\documentclass[11pt]{article}
\usepackage{amsmath}
\usepackage{amssymb}
\usepackage{amsthm}
\usepackage{amsfonts}
\usepackage{graphicx}
\usepackage{pdfpages}
\usepackage{subcaption}
\usepackage{xcolor,colortbl}
\usepackage{multicol}
\usepackage{url}

\definecolor{Gray}{gray}{0.92}
\definecolor{LightCyan}{rgb}{0.88,1,1}

\newcolumntype{a}{>{\columncolor{Gray}}c}
\newcolumntype{b}{>{\columncolor{white}}c}

\usepackage{bbm,epsfig,graphics,epic,color,rotating,color}
\numberwithin{equation}{section}

\def\P{{\mathbb P}}

\def\E{{\mathbb E}}

\usepackage{mathtools}

\DeclarePairedDelimiter\floor{\lfloor}{\rfloor}

\theoremstyle{plain}
\newtheorem{theorem}{Theorem}

\newtheorem{definition}{Definition}
\newtheorem{example}{Example}

\theoremstyle{definition}

\title{The diffusion of opposite opinions in a randomly biased environment}
\author{Manuel Gonz\'alez-Navarrete\thanks{Departamento de Estad\'istica, Universidad del B\'io-B\'io, Chile. e-mail: magonzalez@ubiobio.cl} ,
Rodrigo Lambert \thanks{Faculdade de Matem\'atica, Universidade Federal de Uberl\^andia and Departamento de Matem\'atica, Universidade Federal da Bahia, Brazil. e-mail: rodrigolambert@ufu.br}
\date{}
}

\begin{document}
\maketitle

\begin{abstract}
We propose a model for diffusion of two opposite opinions. Here, the decision to be taken by each individual is a random variable which depends on the tendency of the population, as well on its own trend characteristic. The influence of the population trend can be positive, negative or non-existent in a random form. We prove a phase transition in the behaviour of the proportion of each opinion. Namely, the mean square proportions are linear functions of time in the diffusive case, but are given by a power law in the superdiffusive regime.
\end{abstract}


\section{Introduction}
\label{sec:intro}

The spreading of a disease, rumour, political position or music preference has been a widely studied theme along the years and can be found in the literature of several scientific research areas. Most of the times these studies involve  {some} mathematical model. Since it is reasonable to assume that these processes are inherently probabilistic, then it is suitable to adopt a stochastic approach. In this line, statistical mechanics principles can be considered. As illustrations, we quote a voter-model approach \cite{NySz} and a random-field Ising interaction \cite{Bou}. Essentially, the main feature is the existence of a phase transition in such models.

Following another direction, Bendor et al. \cite{BHW} introduced a model for the diffusion of opinions as a probabilistic choice process. This model has its roots on Festinger's hypothesis \cite{Fe} and is formulated as follows. Two alternatives, $\mathcal{A}$ and $\mathcal{B}$, diffuse through a population (or \emph{decision makers}, following the author's terminology). It is assumed that there is a positive initial number of ``converted" decision makers. In every period, one decision maker makes up his (her) mind about whether to adopt $\mathcal{A}$ or $\mathcal{B}$. The decision taken at the $n$-step is denoted by a Bernoulli variable $X_n$. In this sense, $X_n=1$ means $n$-decision was the innovation $\mathcal{A}$ (and $X_n=0$ means $\mathcal{B}$-decision). Let us denote by $N_n= N_0 + X_1 + \ldots + X_n$ the number of decision makers which adopted $\mathcal{A}$ and $M_n$ those that adopted $\mathcal{B}$, before the $(n+1)$-th decision. Therefore, for $n \geq 1$ the BHW (Bendor-Huberman-Wu) model is defined by the conditional probabilities
\begin{equation}
\label{Ben}
\P\left(X_{n+1} = 1 | (N_n,M_n)\right)=\theta p +(1-\theta)\frac{N_n}{N_n+M_n} \ ,
\end{equation}
where $0< \theta \le 1$ and $P(X_1=1)={N_0}/{(N_0+M_0)}$. When $\theta=0$, we recover the well-known classical P\'olya-urn process (for details, see \cite{Mah}). In another case, if $\theta = 1$, $\{X_i\}_{i \in \mathbb{N}}$ is an independent and identically distributed (i.i.d.) Bernoulli process and then $N_n$ follows a binomial distribution. We remark, as exposed in \cite{BHW}, that each decision is taken by a linear combination of own and social influences, being represented by parameter $p$ and proportion ${N_n}/{(N_n+M_n)}$, respectively.

The role of parameters $p$ and $\theta$ in the mean proportion of opinions $\mathcal{A}$ (namely $\E\left[{N_n}/{(N_n+M_n)}\right]$), is discussed in some details in \cite{BHW}. Particular cases from \eqref{Ben} include behaviours such as herding (see \cite{Ba,BHW2}) and the so-called soft technologies when $p \sim 1/2$ and $\theta$ is close to zero (see \cite{BHW}).

In a different framework, Drezner and Farnum \cite{DF} introduced a generalized binomial distribution. In that model, the independence between Bernoulli trials is replaced by an assumption similar  {to} \eqref{Ben}. By using another parametrization and $P(X_1=1)=p$, the authors showed that, as with the case of the binomial distribution, $\E(X_n)=p$ and $\text{Var}(X_n)=p(1-p)$, for all $n\ge 1$. However, $\text{Var}(N_n)$ is not the same  {as} in the binomial distribution.

We remark that this family of models finds applications in different areas, such as, epidemiology \cite{Ze} or train accidents prevention \cite{LH,LH2}, among others. Moreover, an important attention has been shown in a class of random walks with longe-range correlation \cite{HK,ST}, expressed by conditional probabilities as in \eqref{Ben}. The relation of such processes with P\'olya-type urns was exposed in \cite{BB} (we refer the reader to Ref. \cite{GL} for further details).

Regarding the asymptotic behaviour of $N_n$, Heyde \cite{He} has obtained a law of large numbers that provides a limiting proportion for $N_n$ (and so for $M_n$). Formally, the author proved that, for all  $\theta > 0$, it holds
\begin{equation}
\label{LLN}
\displaystyle\lim_{n\to \infty} \frac{N_n}{N_n+M_n} = p \ \ \text{ a.s. }
\end{equation}
In the same paper it was stated a central limit theorem whenever $\theta \ge 1/2$, and the existence of a phase transition at $\theta = 1/2$. These results are presented in a summarized form in the next theorem. We state this results in terms of the BHW model and present it here without proof.

\begin{theorem}(Heyde, 2004 \cite{He})
\label{PropoA}
Let $\{X_i\}_{i \ge 1}$ defined by conditional probabilities \eqref{Ben} and $P(X_1=1)=p$. The following assertions hold.

\begin{enumerate}
\item [(i)] If $\theta>  \frac{1}{2}$,
 \begin{equation}
\label{TCL}
\frac{N_n - np}{\sqrt{n}} \xrightarrow{d} \mathcal{N}\left(0,\frac{p(1-p)}{2\theta-1}\right),  \ \ \ \text{ as } \  n \to \infty.
\end{equation}

\item [(ii)] If $\theta=  \frac{1}{2}$,
 \begin{equation}
\label{TCLB}
\frac{N_n - np}{\sqrt{n\log n}} \xrightarrow{d} \mathcal{N}(0,p(1-p)),  \ \ \ \text{ as } \  n \to \infty.
\end{equation}

\item [(iii)] If $0<\theta< \frac{1}{2}$,
 \begin{equation}
\label{TCLC}
\frac{N_n - np}{n^{1-\theta}} \xrightarrow{d} W,  \ \ \ \text{ as } \  n \to \infty,
\end{equation}
where $W$ is a proper random variable such that
\begin{equation}
\label{EW}
\E(W)=0, \ \ \E(W^2)=\frac{p(1-p)}{(1-2\theta)\Gamma(2(1-\theta))} \ \text{ and } \ \E(W^3) \neq 0 \ \text{ if } \ p \neq 1/2.
\end{equation}
\end{enumerate}
\end{theorem}

Here the notation $\xrightarrow{d}$ denominates convergence in distribution, and $Z \sim \mathcal{N}(\mu, \sigma^2)$ means that $Z$ is a normal random variable with mean $\mu$ and variance $\sigma^2$.

We remark some features about the results in this theorem. The \emph{diffusive} case $(i)$ is $\sqrt{n}$-scaled and the \emph{critical} one $(ii)$ has $\sqrt{n \log n}$ as scale factor. The \emph{superdiffusive} regime (iii) converges under the $n^{1-\theta}$ scaling. We recall that in the diffusive and critical cases the limiting distribution is a normal random variable, instead of {the superdiffusive regime, which presents a quite different limiting behaviour.}    
Similar results than $(i)-(ii)$ of this theorem was discussed in \cite{BHW}. The rigorous proof of these results can be viewed in \cite{He} and \cite{Dr}.

In the present paper we propose a generalized approach for the spreading of two opposite opinions or preferences. The main contribution is the fact that each decision maker presents a random bias, which can be positive, negative or neutral. In words, it means that any decision maker has a tendency factor, which is inherent to each individual, and works as a latent process that imposes a random ambiance for this diffusion. This particular feature characterizes the so-called \emph{random environment} in which the diffusion occurs.

For this model, we prove a law of large numbers for the proportion of each decision on the population (Theorem \ref{LGN}). Moreover, in the diffusive and critical cases the limiting distributions are showed to be normal random variables and a functional limit theorem can be also obtained (Theorems \ref{CLT} and \ref{continuous}). Instead, the superdiffusive region exhibits non-normal limiting distribution which depends on the initial proportion of opinions (Theorem \ref{CLT} (iii)), then a particular case is studied (Example \ref{example}).

\section{The diffusion model and main results}
\label{se:BSRD}

Imagine a population in which two opposite opinions $\mathcal{A}$ and $\mathcal{B}$ spread. Suppose that this diffusion process starts with $N_0+M_0$ \emph{seeders}, where $N_0$ is the initial number of individuals which has opinion $\mathcal{A}$, and $M_0$ those that have opinion $\mathcal{B}$. Each individual in the population has one of the following three intrinsic characteristics: \emph{trend-follower}, \emph{against-trend} or \emph{indifferent}. This intrinsic characteristic is independently imputed to each individual and follows a ternary random variable, defined as follows.

\begin{definition}\label{trend} For $\alpha, \beta \in (0,1)$ such that $0\leq \alpha+\beta\leq1$, the trend process $\{Y_n\}_{n \in \mathbb{N}}$ is an i.i.d. sequence of ternary random variables (which is also independent of the current proportions of $\mathcal{A}$ and $\mathcal{B}$ decisions on the population) such that, for each $n \ge 1$
$$
\P(Y_{n} = y) = \left\{\begin{array}{cccc}
                     \alpha   & \mbox{if} & y= +1 & \mbox{(Trend-follower),}\\ 
                     \beta   & \mbox{if} & y=-1 & \mbox{(Against-trend),}\\
                       1 - \alpha - \beta & \mbox{if} & y=0 & \mbox{(Indifferent).}
                       
\end{array}
\right.
$$
\end{definition}

We are then ready to introduce the diffusion model with randomly biased environment.
\begin{definition}
\label{our}
The random-trend diffusion model is a stochastic process $\{X_n\}_{n \in \mathbb{N}}$ taking values on $\{0,1\}^{\mathbb{N}}$ defined by the following conditional probabilities
\begin{equation}
\label{Mine}
\P(X_{n+1} = 1 | \mathcal{F}_n, Y_n)=a + bY_n\frac{N_n}{N_n+M_n},
\end{equation}
where $N_n$ (respec. $M_n$) is given by the number of $\mathcal{A}$-decisions (respec. $\mathcal{B}$-decisions) taken until step $n$, $\mathcal{F}_n = \mathcal{F}_n(N_n,M_n)$ is the filtration related to random variables $N_n, M_n$. In view of Definition \ref{trend}, we see that $Y_n$ is independent of $\mathcal{F}_n$, for all $n \ge 1$. Parameters $a$ and $b$ satisfy $0 \le a + b \le 1$, and also the following condition: if $\beta\neq 0$, then $b\leq a$.
\end{definition}

In this sense, the trend process $\{Y_n\}_{n \in \mathbb{N}}$ can be seen as a latent process. In words, it is a process that is not possible to see its input, but that influences the output of each realization of $\{X_n\}_{n \in \mathbb{N}}$. In what follows, we present an example of how this diffusion can be seen. It is a version of an application presented in \cite{GL} and illustrates a potential application for the randomly biased diffusion model.

\begin{example}\label{ex1}
At any given time $n$, a customer can buy a product $(X_n= 1)$ by necessity or by social pressure (history of selling). In the present diffusion model, we consider the conditional probabilities \eqref{Mine}. Let be parameter $a$ the quantity related to necessity and parameter $b$ the quantity related to the social pressure of the product as defined in \cite{BHW,ZH}, associated to $(n+1)$-th customer. 
Then, if $Y_n=1$, the individual is a trend-follower, and the social pressure will raise the probability of the customer buying the product. If $Y_n=-1$, the costumer is against the trend, then the social pressure will decrease the probability of he/she buying the product. Finally, if $Y_n=0$, the costumer behaves without taking the social pressure into account.
\end{example} 

As a particular case, when $\beta=0$ we recover the dependent random walks with memory lapses introduced in \cite{GL}. We say that a $\ell$-sized memory lapse is given by a string $Y_i^{i+ \ell}=0^{\ell}$, in which $Y_i^{i+ \ell}$ is the short notation for $Y_i, Y_{i+1}, \ldots, Y_{i+\ell}$, and $0^{\ell}$ means $\ell$ consecutive zeros.  In a memory lapse situation, the process is occurring independently of the past. When the lapse stops, the process recovers its whole memory. This includes the piece of memory built at the lapse period too (we refer the reader to \cite{GL,GL2} for more details, and to \cite{Be,SB,Ya} to similar ideas).


In what follows we state our first result. It is a strong law of large numbers which provides the limiting proportion of individuals which has taken decisions $\mathcal{A}$ or $\mathcal{B}$, as $n$ diverges.

\begin{theorem}\label{LGN}
Consider the random-trend diffusion model as in Definition \ref{our}. If $0 \le a+b \le 1$ and $b\le a$ we get

\begin{equation}
\label{teo1}
\lim_{n \to \infty}\left(\frac{N_n}{N_n+M_n},\frac{M_n}{N_n+M_n}\right) = \left(\frac{a}{1-b(\alpha-\beta)},\frac{1-a-b(\alpha-\beta)}{1-b(\alpha-\beta)}\right)  \ \ \mbox{a.s.}
\end{equation}
\end{theorem}

In the case $\beta=0$, this theorem holds for all $a,b$ such that $0 \le a+b \le 1$. In the following results this also happens for $\beta=0$.


 {Now, we state a phase transition from diffusive to superdiffusive behaviours. In this sense, the parameter space is essentially divided by a critical line which depends on the values of $b$ and the difference $\alpha -\beta$. As will be shown, each case demands a particular time-scaling.}

\begin{theorem}\label{CLT} Suppose that the hypothesis of Theorem \ref{LGN} are satisfied

\begin{itemize}
\item[(i)] If $b(\alpha-\beta) < \frac{1}{2}$, then
\begin{equation}
\label{teoi}
\frac{1}{\sqrt{n}}\left[(N_n,M_n) - n\left(\frac{a}{1-b(\alpha-\beta)}, \frac{1-a-b(\alpha-\beta)}{1-b(\alpha-\beta)}\right) \right] \xrightarrow{d} \mathcal{N}(0,\Sigma_1) \ ,
\end{equation}
with covariance matrix $\Sigma_1$ given by
\begin{equation}
\label{sigma1}
\Sigma_1 = \dfrac{a(1-a-b(\alpha-\beta))}{(1-b(\alpha-\beta))^2(1-2b(\alpha-\beta))} \left( 
\begin{array}{cc}
1
& 
-1
\\
-1
& 
1
\end{array}
\right) \ .
\end{equation}

\item[(ii)] If $b(\alpha-\beta) = \frac{1}{2}$, then

\begin{equation}
\label{teoii}
\frac{1}{\sqrt{n\log(n)}}\left[(N_n,M_n) - n\left(\frac{a}{1-b(\alpha-\beta)}, \frac{1-a-b(\alpha-\beta)}{1-b(\alpha-\beta)}\right) \right] \xrightarrow{d} \mathcal{N}(0,\Sigma_2) \ ,
\end{equation}

where the covariance matrix $\Sigma_2$ is given by
\begin{equation}
\label{sigma2}
\Sigma_2 = 2a(1-2a)\left( 
\begin{array}{cc}
1
& 
-1
\\
-1
& 
1
\end{array}
\right) \ .
\end{equation}

\item[(iii)]\label{superdiffusive}
Suppose that $ \frac{1}{2}< b(\alpha - \beta) < 1$. Then we have tightness for the sequence of two-dimensional random variables $\{Z_n\}_{n\geq2}$ given by
$$
Z_n := \frac{1}{n^{b(\alpha - \beta)}}\left[{(N_n, M_n) - n\left(\frac{a}{1-b(\alpha-\beta)}, \frac{1-a-b(\alpha-\beta)}{1-b(\alpha-\beta)}\right)}\right] \ .
$$
Furthermore, there exists {a} random bi-dimensional vector $\hat{W} = (\hat{W_1}, \ \hat{W_2})$ such that, as $n$ diverges

$$
Z_n  \longrightarrow \hat{W} \ \ \  \mbox{a.s.}
$$

\end{itemize}
\end{theorem}

 {Note that $(i)$ and $(ii)$ are versions of the central limit theorem. In the case of superdiffusive regime in $(iii)$}, we recall that condition $b(\alpha-\beta)>1/2$ implies that $\alpha$ is a bit larger than  $\beta$. In view of the trend process $\{Y_n\}_{n \in \mathbb{N}}$, one can foretell that the positive trend is expected to be kept along the diffusion. In other words: the statement of this Theorem exhibits an explicit dependence on the initial composition of the model. This result is related to a class of models called large urns, which are usually difficult to obtain explicit limit quantities (for details, see \cite{CPS}).

 {In what follows, we state particular cases to exemplify the previous theorem.}
 {
\begin{example}
\label{example1}
Consider a random-trend model as in Definition \ref{our}, with $\beta=1$ and $a=\frac{1+b}{2}$. Then, we obtain the so-called critical Polya urn introduced in \cite{Ya}. For this model, there is no phase transition and the asymptotic behaviour is characterized by a normal random variable given in item $(i)$.
\end{example}

\begin{example}
\label{example}
Let $\beta=0$ and $\alpha=\frac{1-2a}{b}$. Then, there is a phase transition at $ a = \frac{1}{4}$. In particular, if we start with a single $\mathcal{A}$ opinion, the random vector $\hat{W}$ in item $(iii)$, satisfies, for $a<\frac{1}{4}$

\begin{itemize}
\item[-] $\E(\hat{W})  = \dfrac{1}{2 \Gamma(2(1 -a))} \left( 1 , -1 \right) $ ;
\item[-] $\E(\hat{W}^2) = \dfrac{1}{2 \Gamma( 3-4a)}  \left(1+\dfrac{(1-2a)^2}{1-4a} \right)\left( 1 , -1 \right)$ and
\item[-] $\E(\hat{W}^3) = \dfrac{1}{2\Gamma( 2(2-3a))}  \left(1+\dfrac{(5-2a)(1-2a)^2}{1-4a} \right)\left( 1 , -1 \right) $ \ .
\end{itemize}
\end{example}

This example is related to the elephant random walk introduced in \cite{ST} } and the results hold from the cumulants obtained in Example \textbf{3.27} in \cite{Janson}. In particular, the asymptotic position of the elephant can be obtained by computing $\hat{W}_1 - \hat{W}_2$.


 {Finally, we present a functional central limit theorem for our model. In fact, it is a continuous-time version of Theorem \ref{CLT} (i)-(ii) which deals with a convergence to a two dimensional Brownian motion and provides an explicit limiting covariance matrix. This result is also known in the literature as Donsker's invariance principle with applications arising in estimation, testing, model selection, goodness of fit, regression, and many other common statistical contexts (for details see \cite{DG,He2}).} We recall that this convergence holds on the function space $D[0, \infty )$ of right-continuous and left-bounded (c\`adl\`ag) functions, also known as Skorokhod space.

\begin{theorem} \label{continuous}
 Suppose that the hypothesis of Theorem \ref{LGN} are satisfied.

\begin{itemize}
\item[(i)]  If $b(\alpha-\beta) <  \frac{1}{2}$ then, for $n \to \infty$, in $D[0, \infty)$

\begin{equation}
\label{teo2i}
\frac{1}{\sqrt{n}}\left[(N_{\floor*{tn}},M_{\floor*{tn}}) - tn\left(\frac{a}{1-b(\alpha-\beta)}, \frac{1-a-b(\alpha-\beta)}{1-b(\alpha-\beta)}\right) \right] \xrightarrow{d} W_t,
\end{equation}
where $W_t$ is a continuous bivariate Gaussian process with $W_0=(0,0)$, $\E(W_t)=(0,0)$ and, for $0 < s \le t$,

\begin{equation}
\label{corri}
\E(W_s W_t^T) = s\left(\frac{t}{s}\right)^{b(\alpha - \beta)} \dfrac{a(1-a-b(\alpha-\beta))}{(1-b(\alpha-\beta))^2(1-2b(\alpha-\beta))} \left( 
\begin{array}{cc}
1
& 
-1
\\
-1
& 
1
\end{array}
\right) \ .
\end{equation}

\item[(ii)]  If $b(\alpha-\beta) =  \frac{1}{2}$ then, for $n \to \infty$, in $D[0, \infty)$

\begin{equation}
\label{teo2ii}
\frac{1}{\sqrt{n^t\log(n)}}\left[(N_{\floor*{n^t}},M_{\floor*{n^t}}) - n^t\left(2a, 1-2a\right) \right] \xrightarrow{d} W_t,
\end{equation}
where $W_t$ is a continuous bivariate Gaussian process with $W_0=(0,0)$, $\E(W_t)=(0,0)$ and, \break for $0 < s \le t$,

\begin{equation}
\label{corrii}
\E(W_s W_t^T) = 2sa(1-2a)\left(\begin{array}{cc}
1 & -1 \\
-1 & 1 
\end{array}
\right) \ .
\end{equation}
\end{itemize}
\end{theorem}

\section{Proofs}\label{proofs}

This section is devoted to provide the paths and explicit proofs for the statements given in the paper. We start by constructing the so-called \emph{random replacement matrix}, which was defined in \cite{Janson} as a generalisation of the \emph{replacement matrix} for general P\'olya urn models (see \cite{Mah}, for instance). The eigendecomposition of the random replacement matrix will then provide the main ingredients for the proofs.

\subsection{Building the random reinforcement matrix}
Now we expose the relation between the random-trend diffusion model and a generalized P\'olya urn with random reinforcement matrix. This construction will be done in order to bring the random-trend problem to an urn model and therefore apply the well-known results for generalized P\'olya urns provided by \cite{Janson}.

Firstly we provide an explicit relationship between the distribution of $(N_n,M_n)_{n \in \mathbb{N}}$ and the distribution of the balls in a two-color urn.  
Then, we construct the \emph{random replacement matrix} for the related generalized P\'olya urn. We recall that here and hereafter we will follow the notation given in \cite{Janson}. In this sense, consider the following two column replacement vectors $\xi_{1}=(\xi_{11}, \xi_{12})^T$ (red) and $\xi_{2}=(\xi_{11}, \xi_{12})^T$ (blue), with $\xi_i \in \{(0,1)^T,(1,0)^T\}$ (a single ball is replaced at each time), and the random replacement matrix given by $M = (\xi_1;\xi_2)$. The urn dynamics evolves as follows. If we choose replacement vector $\xi_1$ to reinforce the urn, it means that we will replace $\xi_{11}$ red balls and $\xi_{12}$ blue balls. Otherwise  we choose vector $\xi_2$, and then replace $\xi_{21}$ red, and $\xi_{22}$ blue balls.
At this time, we will be interested in the so-called \emph{mean replacement matrix}, which is given by (with some abuse of notation)

$$
A = \E(M) = \left(
\begin{array}{cc}
\E(\xi_{11}) & \E(\xi_{21}) \\
\E(\xi_{12}) & \E(\xi_{22})
\end{array}
\right) \ .
$$


Note that $\P(\xi_{ij}=1) = \E(\xi_{ij})$, for all $i,j$. If at time $n$ we get $r$ red balls ($R_n=r$), then the probabilities of replacing a red ball (success) or a blue ball (failure) at time $n+1$ (conditioned on a proportion $r/T_n$ of red balls, with $T_n$ the total number of balls at time $n$) are, respectively, given by 
\begin{equation}\label{red}
\begin{array}{ll}
\P(R_{n+1}=r+1 | R_n=r) &= \P(\xi_{11}=1)\dfrac{r}{T_n} + \P(\xi_{21}=1)\left(1-\dfrac{r}{T_n}\right) \\
&= \E[\xi_{21}] + \left(\E[\xi_{11}]-\E[\xi_{21}]\right)\dfrac{r}{T_n},
\end{array}
\end{equation}
and
\begin{equation}\label{blue}
\begin{array}{ll}
 \P(R_{n+1}=r| R_n=r) & = 1 - \P(R_{n+1}=r+1 | R_n=r)\\
                      & = \E[\xi_{22}] + \left(\E[\xi_{12}]-\E[\xi_{22}]\right)\dfrac{r}{T_n} \ .
\end{array}
\end{equation}

In the case of the random-trend diffusion model in Definition \ref{our}, note that, since $Y_n$ and $(N_n, M_n)$ are independent we get that the \emph{mean probability of a success} is given by
\begin{equation}
\begin{array}{lll}
\label{suc}
\P(X_{n+1}=1 | \mathcal{F}_n) & = & \displaystyle\sum_{y \in \{-1,0,1\}} \P(X_{n+1}=1 | \mathcal{F}_n,Y_n=y) \P(Y_n=y) \\[0.5cm]
 & = & a + b(\alpha - \beta) \frac{N_n}{N_n+M_n}. 
\end{array}
\end{equation}
We also get that the \emph{mean probability of a failure} is given by
\begin{equation}\label{fai}
\P(X_{n+1}=0 | \mathcal{F}_n) = 1-\P(X_{n+1}=1 | \mathcal{F}_n) =1-a - b(\alpha - \beta) \frac{N_n}{N_n+M_n}.
\end{equation}
Now we remark the following: in \eqref{suc} and \eqref{fai}, the term $N_n/(N_n+M_n)$ plays the rule of the term $r/T_n$ in \eqref{red} and \eqref{blue}. This correspondence allows us to obtain the entries of the matrix $A$ by simply solving the following system:

$$
\left\{
\begin{array}{l}
\E[\xi_{21}] = a, \\ 
\E[\xi_{22}] = 1-a, \\
\E[\xi_{11}] - \E[\xi_{21}] = b(\alpha-\beta),\\
\E[\xi_{12}] - \E[\xi_{22}] = - b(\alpha-\beta).
\end{array}
\right.
$$

Then we obtain the following matrix

$$
A=  \left(
\begin{array}{cc}
a + b (\alpha - \beta)  & a \\
 1-a-b (  \alpha - \beta)   & 1- a
\end{array}
\right) \ .
$$
The eigenvalues of $A$ are given (in decreasing order) by

$$
\lambda_1= 1  \text{ and } \lambda_2=b(\alpha-\beta) \ .
$$ 

Moreover, the correspondent right and left eigenvectors (which normalizations provided by equations (2.2) and (2.3) of \cite{Janson}) are given by
$$
v_1 = \frac{1}{1-b(\alpha-\beta)}\left(\begin{array}{cc} a \\ 1-a-b(\alpha-\beta) \end{array}\right) \ \ ;  \ \ v_2 = \frac{1}{1-b(\alpha-\beta)} \left(\begin{array}{cc} 1 \\ -1 \end{array}\right),
$$

$$
u_1= (1 \ , \ 1) \ \ \text{ and } \ \ u_2=\left({1-a-b(\alpha-\beta)}\ , \ -a\right).
$$

Now let's check conditions (A1)-(A6) in Section 2 of \cite{Janson}. Since each entry of $A$ is the expected value of a Bernoulli random variable in $\{0,1\}$, we get (A1)-(A2). The facts that $\lambda_1>0$ and $\lambda_2<\lambda_1$ provide (A3)-(A4). Since both colors belong to a dominating class and the urn starts with a positive number of balls, (A5)-(A6) are then satisfied. Moreover, since the urn is tenable, and since it is impossible to exclude balls in the present dynamics, we get that the \emph{essentially non-extinction} condition is also satisfied.

\subsection{Proof of Theorem \ref{LGN}}

\begin{proof}
Since $\lambda_2 = b(\alpha-\beta) < 1 = \lambda_1$ we apply Theorem \textbf{3.21} of \cite{Janson}. It says that if $R_n$ (respect. $B_n$) is the number of red (respect. blue) balls at time $n$, then for $U_n := (R_n, B_n)$ we get
$$
\lim_{n \to \infty}\frac{U_n}{n} = \lambda_1 v_1 \ \ \mbox {a.s.}
$$
 
Since $M_n+N_n = N_0+M_0 + n$, we get that $M_n + N_n \approx n$ when $n$ is big. This completes the proof.
\end{proof}

\subsection{Proof of Theorem \ref{CLT}}

\begin{proof}{(i)}

Theorem \textbf{3.22} of \cite{Janson} to our case says that the limiting covariance matrix is given by
\begin{equation}\label{integral}
\Sigma_1 = \int_{0}^{\infty}\psi_A(s) B \psi_A(s)^T e^{- s}ds -  v_1 v_1^T \ ,
\end{equation}
which in turn depends on the following quantities
\begin{itemize}
\item[(I)] $\psi_A(s) = e^{sA} - v_1 (1 \ 1)\int_{0}^{s}e^{tA} dt$ .

\item[(II)] $B=\sum_{i=1}^2 v_{1i}B_i$ , with $B_i = \E[\xi_i \xi_i^T]$ .
\end{itemize} 

We start by computing term (II). By definition, $B = v_{11}B_1 + v_{12}B_2$, where (again with some abuse of notation)
$$
B_i =  \E[\xi_i \xi_i^T] = \left(
\begin{array}{cc}
\E[\xi_{i1}^2] & \E[\xi_{i1}\xi_{i2}] \\
\E[\xi_{i2}\xi_{i1}] & \E[\xi_{i2}^2]
\end{array}
\right)
=
\left(
\begin{array}{cc}
\E[\xi_{i1}] & 0 \\
0 & \E[\xi_{i2}]
\end{array}
\right) \ ,
$$
where $\xi_i^T$ is the transpose of $\xi_i$, and the last equality comes from the following facts: for all $i\in{1,2}$ $\xi_{ij}\in \{0,1\}$ and the replacement rule does not allow $\xi_{i1}$ and $\xi_{i2}$ being $1$ at the same time. By a direct computation we conclude that

\begin{equation}\label{Bmatrix}
B= \left(\begin{array}{cc}  \frac{a}{1-b(\alpha-\beta)} & 0 \\
0 & \frac{1-a-b(\alpha-\beta)}{1-b(\alpha-\beta)} \end{array}\right).
\end{equation}

For term $(I)$, we use diagonalization $A = V \Lambda V^{-1}$, where $V$ is a diagonal matrix and then $e^{A} = V e^{\Lambda} V^{-1}$. In this case,

$$
V=  \left(
\begin{array}{cc}
1 &  \frac{a}{1-b(\alpha-\beta)} \\
1 &  \frac{1-a-b(\alpha-\beta)}{1-b(\alpha-\beta)}
\end{array}
\right) \ \ \text{and} \ \  \Lambda=  \left(
\begin{array}{cc}
1 & 0 \\
0  & b (  \alpha - \beta)
\end{array}
\right),
$$
allowing us to obtain $\psi_A(s)$ directly.

Since the condition $b(\alpha-\beta)<1/2$ provides convergence for the integral in \eqref{integral}, we use a mathematical software to compute this and then the proof is finished.
%
%

(ii)
Since $A$ is diagonalizable, Theorem \textbf{3.23} of \cite{Janson} to our case says that the limiting covariance matrix is given by

\begin{equation}\label{critical}
\Sigma_2 = (I-T_1) P_{\frac{1}{2}}B P_{\frac{1}{2}}^*(I-T_1)^T \ ,
\end{equation}
where $P_{\frac{1}{2}}^* = P_{\frac{1}{2}}^T$ is the Hermitian conjugate of $P_{\frac{1}{2}}$ (and the equality holds since it is a real-entries matrix). Moreover $I$ is the two-dimensional identity matrix and $T_1$ is given by
\begin{equation}\label{T1}
T_1 = v_1 . (1 \ 1) 
= \frac{1}{1-b(\alpha-\beta)}\left(
\begin{array}{cc}
a &  a \\ 
1-a-b(\alpha-\beta) & 1-a-b(\alpha-\beta)
\end{array}
\right) \ .
\end{equation}

We note that the above equality was obtained by recalling that, in view of \cite{Janson}, since $det(A^T) = 1/2 \neq 0$, $A^T$ has full rank, and then the \emph{action vector} (which in our case is given by $\overrightarrow{a}=(1 \ 1)$) satisfies $\overrightarrow{a} \in Im(A^T)$. Moreover  $\overrightarrow{b}=(1 \ 1)$ satisfies the condition $\overrightarrow{a} = A^T \overrightarrow{b}$. 
Furthermore the fact that $\Lambda_{II} = \{\lambda_2\}  = \{b(\alpha-\beta)\} =  \{1/2\}$ leads to

\begin{equation}\label{P1/2}
P_{\frac{1}{2}} = v_2 u_2 =  2\left(
\begin{array}{cc}
\frac{1}{2}-a & a \\
a - \frac{1}{2} & -a 
\end{array}
\right) \ .
\end{equation}

Since \eqref{Bmatrix}, \eqref{T1} and \eqref{P1/2} provide all the ingredients, we finally compute \eqref{critical}.

 {
Finally, for item (iii) we just need to notice that the eigenspace $E_{\lambda_2}$ is the set of vectors on the direction of $(1,-1)$ and then apply Theorem \textbf{3.24} of \cite{Janson}. Theorem \textbf{3.26} of \cite{Janson} provides the moment generator function for the Example \ref{example}.
}
\end{proof}

\subsection{Proof of Theorem \ref{continuous}}
\begin{proof}
(i) (The diffusive case)

The convergence statement is directly provided by Theorem $\textbf{3.31}-(i)$ of \cite{Janson}, with limiting covariance matrix given by
$$
\E(W_s W_t^T)=\int_{-\log s}^{\infty} \psi_A(x+\log s) B \psi_A(x+\log t) e^{-x} dx - s v_1 v_1^T,
$$
which is equivalent to obtain (see Remark \textbf{5.7} from \cite{Janson})
\begin{equation}
\E(W_s W_t^T) = s \Sigma_1 e^{\log(t/s) A^T}.
\end{equation}

Now, as done in the proof of Theorem \ref{CLT}, we obtain $A^T = V \Lambda V^{-1}$, and then $e^{A^T} = V e^{\Lambda} V^{-1}$. In this case, they are given by

$$
V=  \left(
\begin{array}{cc}
1 & \frac{a+b(\alpha - \beta) -1}{a} \\
1 & 1
\end{array}
\right), 
\ \ \ \Lambda=  \left(
\begin{array}{cc}
1 & 0 \\
0  & b (  \alpha - \beta)
\end{array}
\right) \ \ \text{and} \ \ 
V^{-1}=  \left(
\begin{array}{cc}
\frac{a}{1-b (  \alpha - \beta)} & \frac{a+b(\alpha - \beta) -1}{b (  \alpha - \beta)-1} \\
\frac{a}{b (  \alpha - \beta)-1} & \frac{a}{1-b (  \alpha - \beta)}
\end{array}
\right).
$$

A direct calculation completes the proof.

(ii) (The critical case)

Once again the convergence statement is directly provided by Theorem $\textbf{3.31}-(ii)$ of \cite{Janson}, with $d:=d_{\lambda_2}=0$, since $A$ is diagonalizable. Now let us compute the limiting covariance matrix. To do this, notice that $\tilde{c}(d,s,t) = s$ and then
$$
\E(W_s W_t^T) = s(I - 2v_1 \mathbf{1}  P_{\lambda_2}) \Sigma_{II} (I - 2 P_{\lambda_2}^T \mathbf{1}^T v_1^T).
$$
Actually we use \eqref{P1/2}. Note that,

$$
\mathbf{1}  P_{\frac{1}{2}} = (0 \ 0) \ ; \ P_{\frac{1}{2}}^T \mathbf{1}^T = (0 \ 0)^T.
$$
Then $\E(W_s W_t^T) =  P_{\frac{1}{2}} B P_{\frac{1}{2}}^T$. Concluding the proof, we use \eqref{Bmatrix} and \eqref{P1/2} to obtain the following limiting covariance matrix
$$
\E(W_s W_t^T) = 2sa(1-2a)\left(\begin{array}{cc}
1 & -1 \\
-1 & 1 
\end{array}
\right).
$$
\end{proof}

\subsection*{Acknowledgements}
The authors thank Durval Tonon and Alejandra Rada-Mora for helpful comments and discussions. MGN is supported by Fondo Especial DIUBB 1901083-RS and FAPEI from Universidad del B\'io-B\'io. RL is partially supported by FAPESP (grant 2014/19805-1), and CNPQ (grants 406324/2017-4, 439422/2018-3, and 435470/2018-3). This work is part of the Universal FAPEMIG project ``Din\^amica de recorr\^encia para shifts aleat\'orios e processos com lapsos de mem\'oria" (process APQ-00987-18).

\end{document}